\documentclass[10pt,a4paper]{article}
\usepackage{makeidx}
\makeindex            % damit eine Indexdatei angelegt wird

\usepackage{amsmath}  % allgemeine Mathe-Erweiterungen
\usepackage{amssymb}  % Symbole und Schriftarten
\usepackage{amsthm}   % erweiterte Theorem-Umgebungen
\usepackage{amsfonts}
\newcommand{\N}{\mathbb{N}}

\usepackage{mathrsfs}  % gibt den Befehl "\mathscr{}" f\"{u}r sch\"{o}ne
\newcommand{\Ann}[0]{\operatorname{Ann}}
\newcommand{\cohdim}[0]{\operatorname{cd}}
\newcommand{\Hom}[0]{\operatorname{Hom}}

\newtheorem{satz}{Theorem}[section]
\newtheorem{lemma}[satz]{Lemma}

\newtheorem{remark}[satz]{Remark}
\newtheorem{remarks}[satz]{Remarks}

\newtheorem{question}[satz]{Question}

\author{M. Hellus}
\title{A note on the vanishing of certain local cohomology modules}
\date{}   % Standard: Datum der Kompilierung ("\today").

\begin{document}

\maketitle

\begin{abstract}

For a finite module $M$ over a local, equicharacteristic ring $(R,m)$, we show that the well-known formula $\cohdim(m,M)=\dim M$ becomes trivial if ones uses Matlis duals of local cohomology modules together with spectral sequences. We also prove a new, ring-theoretic vanishing criterion for local cohomology modules.

\end{abstract}
\section{Introduction}
Let $R$ be a noetherian ring, $I$ an ideal of $R$ and $M$ an
$R$-module; one denotes the $n$-th local cohomology module of $M$
with respect to $I$ by $H^n_I(M)$ and the cohomological dimension of $I$ on $M$ by
\[ \cohdim(I,M):=\sup \{ l|H^l_I(M)\neq 0\} .\]
From now on assume that $(R,m)$ is local and $M$ is finitely generated. Grothendieck's Vanishing Theorem (VT) says that $\cohdim (I,M)\leq \dim M$ and Grothendieck's Non-Vanishing Theorem (NVT) says $H^{\dim M}_m(M)\neq 0$. Both are well-known theorems with various proofs, see e.~g. \cite[Theorem 6.1.2]{brodmann98}, \cite [Theorem 2.7]{hartshorne77} (a version for sheaves) for VT and \cite[Theorem 6.1.4]{brodmann98}, \cite[Theorem 7.3.2]{brodmann98} for NVT. The case $I=m$ of VT and NVT {\it together} say that the cohomological dimension is the Krull dimension:
\[ \tag{$*$}\cohdim (m,M)=\dim M.\]

The first aim of this paper is to show that, using Matlis duals of local cohomology modules, formula $(*)$ become almost trivial once one knows:

\begin{enumerate}

\item[$(A)$] \label{lcm_koszul}The fact that local cohomology can be written as
the direct limit of Koszul cohomologies; it is an easy exercise to check that immediate consequences of this are

\begin{enumerate}

\item[$(A_1)$] the base-change formula
$_RH^i_{IS}(N)=H^i_I(_RN)$ ($S/R$ a noetherian algebra, $N$ an
$S$-module, $I$ an ideal of $R$ and $i\in \N$)

\item[$(A_2)$] the formula
\[ H^j_{(X_1,\dots ,X_i)}(k[[X_1,\dots ,X_i]])=\begin{cases} 0, & \mbox{if  }j>i\\ E_{k[[X_1,\dots ,X_i]]}(k)=k[X_1^{-1},\dots
,X_i^{-1}], & \mbox{if }j=i\end{cases}
\]
($k$ a field, $X_1,\dots ,X_i$ indeterminates)

\item[$(A_3)$] the fact that each local cohomology functor of the form $H^j_{(x_1,\ldots ,x_i)R}$ is zero for $j>i$; in particular, $H^i_{(x_1,\ldots ,x_i)}R$ is right exact.

\end{enumerate}

\item[$(B)$] Some Matlis duality theory and some spectral sequence theory. Both serve as {\it technical tools}.

\end{enumerate}

Our method works {\it only in the equicharacteristic case}.

The second aim is to prove theorem \ref{two}, which is a new (sufficient) criterion for the vanishing of local cohomology modules, which is of a ring-theoretic nature; the idea which is used in its proof is, to the best of our knowledge, completely new in this context.

\section{(Non-)Vanishing Theorem}

Everything in this paper is based on the following easy

\begin{lemma}

\label{one} Let $(R,m)$ be a noetherian local complete ring
containing a field $k$, $M$ an $R$-module and $x_1,\dots ,x_i\in R$.
Then
\[ H^i_{\underline xR}(M)\neq 0\iff \dim (R_0)=i\text{ and }\Hom
_{R_0}(M,R_0)\neq 0\] where $R_0:=k[[x_1,\dots ,x_i]]$ as a subring
of $R$ and $\underline x:=x_1,\dots ,x_i$.

\end{lemma}

{\it Proof. } $\Rightarrow $: Assume $\dim (R_0)<i$. Write
$R_0=k[[X_1,\dots ,X_i]]/I$ where $X_1,\dots ,X_i$ are
indeterminates and $I$ is a non-zero ideal of $k[[X_1,\dots
,X_i]]=:S$. Then
\[ H^i_{\underline xR_0}(R_0)\buildrel (A_1),(A_3)\over=H^i_{\underline XS}(S)\otimes
_S(S/I)=0\] as every $0\neq f\in I$ operates injectively on $S$ and
hence ($(B)$) surjectively on $H^i_{\underline XS}(S)$($\buildrel (A_2)\over\cong E_S(k)$). In
particular, \[H^i_{\underline xR}(M)\buildrel{(A_3)}\over=M\otimes _{R_0}H^i_{\underline
xR_0}(R_0)=0,\] contradiction. Therefore, $\dim (R_0)=i$, $R_0\cong
k[[X_1,\dots ,X_i]]$ with indeterminates $X_1,\dots ,X_i$ and one
has
\begin{eqnarray*} 0&\buildrel (B)\over \neq &\Hom _{R_0}(H^i_{\underline
xR}(M),E_{R_0}(k))\\ &\buildrel (A_3)\over=&\Hom _{R_0}(M\otimes _{R_0}H^i_{\underline
xR_0}(R_0),E_{R_0}(k))\\ &=&\Hom _{R_0}(M,\Hom
_{R_0}(H^i_{\underline xR_0}(R_0),E_{R_0}(k)))\\ &\buildrel (A_2),(B)\over=&\Hom
_{R_0}(M,R_0)\\
\end{eqnarray*}

$\Leftarrow $: Again, $R_0\cong k[[X_1,\dots ,X_i]]$ with
indeterminates $X_1,\dots ,X_i$; now,
\[ 0\neq \Hom _{R_0}(M,R_0)=\Hom _{R_0}(H^i_{\underline
xR}(M),E_{R_0}(k))\] follows like above.\hfill $\square $

\begin{satz}

\label{cd_eq_dim}(i) If $R$ is a noetherian ring containing a field,
$\underline x=x_1,\dots ,x_i\in R$ and $M$ is an $R$-module (not
necessarily finitely generated) such that $\dim _R(M)<i$, then
$H^i_{\underline xR}(M)=0$.

(ii) If $(R,m)$ is a noetherian local ring containing a field and
$\underline x=x_1,\dots ,x_i$ is part of a system of parameters of a
finitely generated $R$-module $M$ then $H^i_{\underline xR}(M)\neq
0$; in particular, $H^{\dim _R(M)}_m(M)\neq 0$.

(iii) If $(R,m)$ is a noetherian local ring containing a field and
$M$ is a finitely generated $R$-module then $\cohdim (m,M)=\dim
_R(M)$.

\end{satz}

{\it Proof. }(i) By localizing and completing we may assume that $R$
is local and complete. Set $R_0:=k[[x_1,\dots ,x_i]]$ as a subring
of $R$ like in lemma \ref{one}; we may assume that $\dim (R_0)=i$,
i.~e. $R_0\cong k[[X_1,\dots ,X_i]]$, where $X_1,\dots ,X_i$ are
indeterminates. Due to dimension reasons it is clear that $\Hom
_{R_0}(M,R_0)=0$ and the claim follows from lemma \ref{one}.

(ii) We may assume that $R$ is complete ($\hat R/R$ is faithfully flat); by base-change, we may replace $R$
by $R/\Ann _R(M)$; set $d:=\dim (R)$. We choose $x_{i+1},\dots
,x_d\in R$ such that $x_1,\dots ,x_d$ is a system of parameters of
$M$. Then $R_0:=k[[x_1,\dots ,x_d]]\subseteq R$ is a regular
$d$-dimensional subring of $R$ and, because $M$ is module-finite
over $R_0$, $\Hom _{R_0}(M,R_0)\neq 0$; lemma \ref{one} implies
$H^d_{(x_1,\dots ,x_d)R}(M)\neq 0$. Now a formal spectral sequence
argument (namely for the spectral sequence of composed functors $E_2^{p,q}=H^p_{(x_{i+1},\ldots ,x_d)R}(H^q_{(x_1,\ldots ,x_i)R}(M))\Rightarrow H^{p+q}_{(x_1,\ldots ,x_d)R}(M);$ note that $H^p_{(x_{i+1},\ldots ,x_d)R}=0$ for each $p>d-i$ and that $H^q_{(x_1,\ldots ,x_i)R}=0$ for each $q>i$, by $(A_3)$ ) shows
\[ 0\neq H^d_{(x_1,\dots ,x_d)R}(M)=H^{d-i}_{(x_{i+1},\dots
,x_d)R}(H^i_{(x_1,\dots ,x_i)R}(M))\] (iii) Follows from (i) and
(ii).\hfill $\square $

\section{A Ring-theoretic Vanishing Criterion}

\label{suff_crit}

\begin{satz}

\label{two} Let $(R,m)$ be a noetherian local complete domain
containing a field and $\underline x=x_1,\dots ,x_i$ a sequence in
$R$. Then the implication
\[ H^i_{\underline xR}(R)\neq 0\Rightarrow \dim (R_0)=i\text{ and }R\cap
Q(R_0)=R_0\] holds, where $R_0:=k[[x_1,\dots ,x_i]]\subseteq R$,
$Q(R_0)$ denotes the quotient field of $R_0$ and the intersection is
taken inside $Q(R)$.

\end{satz}

{\it Proof. }By lemma \ref{one}, $R_0\cong k[[X_1,\dots ,X_i]]$ ,
$X_1,\dots ,X_i$ indeterminates, $\dim (R_0)=i$.

Let $r\in R, r_0\in R_0$ such that $r_0\cdot r\in R_0$. We have to
show that $r\in R_0$: by lemma \ref{one}, $\Hom _{R_0}(R,R_0)\neq
0$ and so we can choose $\varphi \in \Hom _{R_0}(R,R_0)$ such that
$\varphi (1_R)\neq 0$ (namely by composing a $\varphi ^\prime \in
\Hom _{R_0}(R,R_0)$ that has $\varphi (r^\prime )\neq 0$ (for some
$r^\prime \in R$) with the multiplication map $R\buildrel r^\prime
\over \to R$). Set $r_0^\prime :=r_0r$. One has
\[ r_0\varphi (r)=\varphi (r_0^\prime )=r_0^\prime \varphi (1_R)\]
and then
\[ \varphi (1_R)r=\varphi (1_R){r_0^\prime \over r_0}=\varphi (r)\in
R_0\] On the other hand, we have \[ r_0^{\prime 2}=r_0^2r^2\] and
thus \[ r_0^2\varphi (r^2)=r_0^{\prime 2}\varphi (1_R)\] and \[
\varphi (1_R)r^2=\varphi (1_R){r_0^{\prime 2}\over r_0^2}=\varphi
(r^2)\in R_0\ \ .\] Continuing in the same way, one sees that, for
every $l\geq 1$, one has \[ \varphi (1_R)r^l\in R_0\ \ .\] But this
implies that the $R_0$-module \[ \varphi (1_R)\cdot <1,r,r^2,\dots
>_{R_0}\] is finitely generated ($<1,r,r^2,\dots >_{R_0}$ stands for
the $R_0$-submodule of $R$ generated by $1,r,r^2,\dots $). But, as
$R$ is a domain, \[ <1,r,r^2,\dots >_{R_0}\] is then finitely
generated, too, i.~e. $r$ is necessarily contained in $R_0$.\hfill
$\square $

\begin{remarks}

\label{eine_remarks}(i) $H^i_{\underline xR}(R)\neq 0$ (and thus
$R\cap Q(R_0)=R_0$) are clear if $\underline x$ is an $R$-regular
sequence; but the condition $\underline x$ being a regular sequence
is not necessary as the following example shows:
$H^2_{(y_1y_2,y_1y_3)}(k[[y_1,y_2,y_3]])$ is non-zero (and thus
$R\cap Q(R_0)=R_0$) though $y_1y_2,y_1y_3$ is not a regular sequence
($k$ a field, $y_1,y_2,y_3$ indeterminates).

(ii) In the situation of theorem \ref{two} without the assumption
$H^i_{\underline xR}(R)\neq 0$ the condition $R\cap Q(R_0)=R_0$ does
not hold in general: e.~g. for $R_0=k[[y_1y_2,y_1y_2^2]]\subseteq
k[[y_1,y_2]]=R$ ($k$ a field, $y_1,y_2$ indeterminates) one has
$y_2\in (R\cap Q(R_0))\setminus R_0$.
\end{remarks}

\begin{remark}

\label{reverse} If $R$ is regular, the implication from theorem
\ref{two} is an equivalence for $i=1$; while this is easy to see, the case $i=2$ seems already
unclear.

\end{remark}

\begin{question}

\label{quest} Under what conditions can the implication from theorem
\ref{two} be reversed?
\end{question}

\end{document}